\documentclass[journal,twoside,web]{ieeecolor}
\usepackage{lcsys}

\newcommand\copyrighttext{%
\footnotesize \textcopyright 2022 IEEE. Personal use of this material is permitted. Permission from IEEE must be obtained for all other uses, in any current or future media, including reprinting/republishing this material for advertising or promotional purposes, creating new collective works, for resale or redistribution to servers or lists, or reuse of any copyrighted component of this work in other works.}
\newcommand\copyrightnotice{%
\begin{tikzpicture}[remember picture,overlay]
	\node[anchor=south,yshift=5pt] at (current page.south) {\fbox{\parbox{\dimexpr\textwidth-\fboxsep-\fboxrule\relax}{\copyrighttext}}};
\end{tikzpicture}%
}

\usepackage{cite}
\usepackage{amsmath,amssymb,mathtools,amsfonts}
\usepackage{graphicx}
\usepackage{textcomp}

\usepackage{xcolor}
\usepackage{listings}
\definecolor{codegreen}{rgb}{0,0.6,0}
\definecolor{codegray}{rgb}{0.5,0.5,0.5}
\definecolor{codepurple}{rgb}{0.58,0,0.82}
\definecolor{backcolour}{rgb}{0.95,0.95,0.92}
\lstdefinestyle{mystyle}{
	basicstyle=\ttfamily\footnotesize,
	breakatwhitespace=false,         
	breaklines=true,                 
	captionpos=b,                    
	keepspaces=true,                 
	numbers=left,                    
	numbersep=5pt,                  
	showspaces=false,                
	showstringspaces=false,
	showtabs=false,                  
	tabsize=2,
	xleftmargin=0.025\textwidth,
	xrightmargin=0.025\textwidth
}
\usepackage{hyperref}
\usepackage{multirow}
\usepackage{color}
\usepackage{psfrag}
\usepackage{subfigure}

\usepackage{empheq}
\usepackage{tabularx}
\usepackage{bm} 
\usepackage{url}
\usepackage{psfrag}
\usepackage{float}
\usepackage{comment}

\usepackage{alphalph}

\usepackage{tikz}
\usepackage{pgfplots}
\usetikzlibrary{arrows.meta}
\newlength\fwidth
\newlength\fheight
\usetikzlibrary{arrows}
\newcommand{\R}{{\mathbb{R}}}
\newcommand{\dd}{{\mathrm{d}}}

\newcommand{\Nsys}{{n_{f}}} 
\newcommand{\Nstg}{{n_\mathrm{s}}} 
\newcommand{\NFE}{{N_{\mathrm{FE}}}} 
\newcommand{\Nctrl}{{N}_{\mathrm{stg}}} 

\newcommand{\I}{{\mathcal{I}}}
\newcommand{\irk}{{\mathrm{rk}}}
\newcommand{\lp}{{\mathrm{LP}}}
\newcommand{\fesd}{{\mathrm{fesd}}}
\newcommand{\NOSNOC}{\texttt{NOSNOC}}

\usepackage{amsthm}

\begin{document}

\title{NOSNOC: A Software Package for Numerical Optimal Control of Nonsmooth Systems }
\author{Armin Nurkanovi\'c and Moritz Diehl
\thanks{This research was supported by the DFG via Research Unit FOR 2401 and project 424107692 and by the EU via ELO-X 953348.	}
\thanks{Armin Nurkanovi\'c is with the Department of Microsystems Engineering (IMTEK), University of Freiburg, Germany,	
		Moritz Diehl is with the
		Department of Microsystems Engineering (IMTEK) and Department of Mathematics, University of Freiburg, Germany,  \texttt{\{armin.nurkanovic,moritz.diehl\}@imtek.uni-freiburg.de}}
}
\maketitle
\thispagestyle{empty} 
\copyrightnotice
\begin{abstract}
This letter introduces the NOnSmooth Numerical Optimal Control (NOSNOC) open-source software package. It is a modular MATLAB tool based on CasADi and IPOPT for numerically solving Optimal Control Problems (OCP) with piecewise smooth systems (PSS).
The tool supports: 1) automatic reformulation of systems with state jumps into PSS (via the time-freezing reformulation \cite{Nurkanovic2021}) and of PSS into computationally more convenient forms, 2) automatic discretization of the OCP via, e.g., the recently introduced Finite Elements with Switch Detection \cite{Nurkanovic2022} which enables high accuracy optimal control and simulation of PSS, 3) solution methods for the resulting discrete-time OCP.
The nonsmooth discrete-time OCP are solved with techniques of continuous optimization in a homotopy procedure, without the use of integer variables. This enables the treatment of a broad class of nonsmooth systems in a unified way. 
Two tutorial examples are given. A benchmark shows that NOSNOC provides both faster and more accurate solutions
than conventional approaches, including mixed-integer formulations.
\end{abstract}
\vspace{-0.25cm}
\section{Introduction}
\IEEEPARstart{N}{onsmooth} and hybrid dynamical systems are a powerful tool to model complex physical and cyber-physical phenomena. Their theory is well established and many good numerical simulation algorithms exist \cite{Brogliato2020}. However, optimal control of nonsmooth systems is yet not wide spread, mainly due to the computational difficulty and lack of software. A notable exception are mixed integer optimization approaches \cite{Bemporad1999b}. However, they become intractable as soon as nonconvexities appear or exact junction times need to be computed. The open-source software package \NOSNOC\ is designed to reduce this gap \cite{Nurkanovic2022c}.
We regard a nonsmooth OCP of the following form:
\begin{subequations} \label{eq:ocp}
	\begin{align}
		\min_{x(\cdot),u(\cdot)} \quad & \int_{0}^{T} f_{\mathrm{q}}(x(t),u(t))\dd t+  f_{\mathrm{T}}(x(T)) \\
		\textrm{s.t.} \quad  x_{0} &= s_0, \\
		\dot{x}(t) &\! = \! f_i(x(t),u(t)),\ \!\mathrm{if}\ \! x\! \in\! R_i,\ \! i \!\in\!\mathcal{I}, t \in [0,\!T], \label{eq:ocp_pss} \\
		0&\geq G_{\mathrm{ineq}}(x(t),u(t)),\ t \in [0,T],\\
		0&\geq G_{\mathrm{T}}(x(T)),
	\end{align}
\end{subequations}
where $f_{\mathrm{q}}: \R^{n_x} \times \R^{n_u} \to \R$ is the stage cost and $f_{\mathrm{T}}:\R^{n_x}\to \R$ is the terminal cost, $s_0\in\R^{n_x}$ is a given parameter. The path and terminal constraints are collected in the functions $G_{\mathrm{ineq}} : \R^{n_x}  \times \R^{n_u} \to \R^{n_{g1}}$ and $G_{\mathrm{T}} : \R^{n_x}  \to \R^{n_{g2}}$, respectively. The ODE \eqref{eq:ocp_pss} is a piecewise smooth system (PSS), where $\mathcal{I} \coloneqq \{ 1,\ldots,\Nsys  \}$. The regions $R_i$ are disjoint, nonempty, connected and open. The functions $f_i(\cdot)$ are smooth on an open neighborhood of $\overline{R_i}$, which denotes the closure of ${R}_i$. 
\color{black}


 The event of $x$ reaching some boundary $\partial R_i$ is called a \textit{switch}. The right-hand side (r.h.s.) of \eqref{eq:ocp_pss} is in general discontinuous in $x$. Several classes of systems with state jumps can be brought into the form of \eqref{eq:ocp_pss} via the time-freezing reformulation \cite{Nurkanovic2021,Nurkanovic2021c,Nurkanovic2022a}. Thus, the focus on PSS enables a unified treatment of many different kinds of nonsmooth systems.

One might wonder why not just to apply standard direct methods and existing software to problem with a smoothed version of the r.h.s. of \eqref{eq:ocp_pss}? The necessity for tailored methods and software follows from two important results from the seminal paper of Stewart and Anitescu \cite{Stewart2010}. First, in standard direct approaches for \eqref{eq:ocp_pss}, the numerical sensitivities are wrong no matter how small the integrator step-size is. This often yields artificial local minima and impairs the optimization progress~\cite{Nurkanovic2020}. Second, smoothing delivers correct sensitivities only if the step-size shrinks faster than the smoothing parameter. Consequently, even for moderate accuracy, many optimization variables are needed.

These two difficulties are overcome by the recently introduced Finite Elements with Switch Detection (FESD) method~\cite{Nurkanovic2022}. In this method, the ODE \eqref{eq:ocp_pss} is transformed into a Dynamic Complementarity System (DCS). FESD relies on  Runge-Kutta (RK) discretizations of the DCS, but the integrator step-sizes are left as degrees of freedom as first proposed by \cite{Baumrucker2009}. Additional constraints ensure implicit and exact switch detection and eliminate spurious degrees of freedom. The discretization yields Mathematical Programs with Complementarity Constraints (MPCC). They are highly degenerate and nonsmooth Nonlinear Programs (NLP) \cite{Scholtes2001,Anitescu2007}, but with suitable reformulations and homotopy procedures they can be solved efficiently using techniques for smooth NLP, without any integer variables.

The MATLAB tool \NOSNOC\ \cite{Nurkanovic2022c} aims to automate the whole tool-chain and to make nonsmooth optimal control problems solvable for non-experts. In particular, it supports:
\begin{itemize}
	 
	\item automatic model reformulation of the PSS \eqref{eq:ocp_pss} into the computationally more suitable DCS. 
	\item time-freezing reformulation for systems with state jumps, reformulations to solve time-optimal control problems both for PSS and systems with state jumps,
	\color{black}
   	\item automatic discretization of the OCP \eqref{eq:ocp} via FESD or RK,
   	\item several algorithms for solving the MPCC with a homotopy approach,    	
   	\item rapid prototyping with different formulations and algorithms for nonsmooth OCP.
\end{itemize}
It builds on the open-source software packages: {CasADi} \cite{Andersson2019} which is a symbolic framework for nonlinear optimization and the NLP solver {IPOPT}\cite{Waechter2006}. Having these packages as a back-end enables good computational performance, despite the fact that all user inputs are provided in {MATLAB}. All steps above can be performed in a couple of lines of code without needing a deep understanding of the numerical methods and implementation details.
In \NOSNOC, the user has only to specify the functions in \eqref{eq:ocp} and the sets $R_i$ via constraint functions $c(x)$, cf. Section \ref{sec:prelimit}. The reformulation, discretization and solution of the nonsmooth OCP is completely automated.
\vspace{-0.4cm}
\paragraph*{Notation}
The complementarity conditions for two vectors  $a,b \in \R^{n}$ read as ${0\leq a \perp b\geq 0}$, where $a \perp b$ means $a^{\top}b =0$. The so-called C-functions $\Phi: \R^{n} \times \R^{n} \to \R^{n}$ have the property $\Phi(a,b) = 0\!\! \iff\!\! {0\leq a \perp b\geq 0}$, e.g., $\Phi(a,b) = \min(a,b)$.
The concatenation of two column vectors $a\in \R^{n_a}$, $b\in \R^{n_b}$ is denoted by $(a,b)\coloneqq[a^\top,b^\top]^\top$, the concatenation of several column vectors is defined in an analogous way. A column vector with all ones is denoted by $e=(1,1,\dots,1) \in \R^n$, its dimensions is clear from the context.
The closure of a set $C$ is denoted by $\overline{C}$, its boundary as $ \partial C$. Given a matrix $M \in \R^{n \times m}$, its $i$-th row is denoted by $M_{i,\bullet}$ and its $j$-th column is denoted by $M_{\bullet,j}$. 
\paragraph*{Outline}
Section \ref{sec:prelimit} describes the reformulation of PSS into DCS. Section \ref{sec:algoritmic_ing} describes the discretization methods in \NOSNOC\, with a focus on FESD. In Section \ref{sec:ocp_discretization}, solution strategies for the discrete-time OCP are discussed. Section \ref{sec:nosnoc} provides two tutorials for the use of \NOSNOC\ and a numerical benchmark. Section \ref{sec:conclusions} outlines some future developments.
\section{Problem Reformulation}\label{sec:prelimit}
System with state jumps do not fit in the form of \eqref{eq:ocp_pss}. However, we use the time-freezing reformulation \cite{Nurkanovic2021,Nurkanovic2022a,Nurkanovic2021c} to automatically reformulate them into the from of \eqref{eq:ocp_pss}. An example is given in Section \ref{sec:time_freezing_example}.
\color{black}

In this section, we detail how to compactly represent the systems \eqref{eq:ocp_pss} and a how to transform them into a Dynamic Complementarity System (DCS) via Stewart's approach \cite{Stewart1990b}.

It is assumed that $\overline{\bigcup\limits_{i\in \mathcal{I}} R_i} = \R^n$ and that $\R^n \setminus \bigcup\limits_{i\in\mathcal{I}} R_i$ is a set of measure zero. Moreover, we assume that $R_i$ are defined via the zero level sets of the components of the smooth function $c:\R^{n_x} \to \R^{n_c}$. We use a sign matrix $S \in \R^{\Nsys \times n_c}$ with non repeating rows for a compact representations
as follows: 
\begin{subequations}\label{eq:standard_sets}
\begin{align}
	S&= \begin{bmatrix}
		1 & 1 & \dots &1 & 1\\
		\vspace{-0.05cm}
		1 & 1 & \dots & 1& -1\\
		\vspace{-0.08cm}
		\vdots &  \vdots & \dots &  \vdots\\
		\vspace{-0.08cm}
		-1 & -1 & \dots &-1 & -1\\
	\end{bmatrix},\\
R_i &= \{ x\in \R^{n_x} \mid \mathrm{diag} (S_{i,\bullet}) c(x)>0\}.
\end{align}
\end{subequations}
 
For example, for the sets $R_1 = \{x\in \R \mid x >0 \}$ and $R_2 = \{x\in \R \mid x <0 \}$, we have $c(x) = x$ and $ S = \begin{bmatrix}
	1 &-1
\end{bmatrix}^\top$.
\color{black}

The dynamics are not defined on $\partial R_i$ and to have a meaningful notion of solution for the PSS \eqref{eq:ocp_pss} we use the Filippov convexification and define the following differential inclusion \cite{Filippov2013}:
\begin{align}\label{eq:FilippovDI_with_multiplers}
	\begin{split}
	\dot{x}  \in    F_{\mathrm{F}}(x,u) = &\Big\{ F(x)\theta  \mid \sum_{i\in \I}\theta_i = 1,
	\ \theta_i \geq 0,
	\ \theta_i = 0 \  \\
	&\mathrm{if} \;  x \notin \overline{R_i}, 
	\forall  i  \in \I \Big\},	
\end{split}
\end{align}
where $\theta = (\theta_1,\ldots,\theta_{\Nsys}) \in \R^{\Nsys}$ and $F(x) \coloneqq [f_1(x),\ldots,f_{\Nsys}(x)] \in \R^{n_x \times \Nsys}$. Note that in the interior of a set $R_i$ we have $F_{\mathrm{F}}(x) = \{f_i(x)\}$ and on the boundary between some regions the resulting vector field is a convex combination of the neighboring vector fields. 
To have a computationally useful representation of the Filippov system~\eqref{eq:FilippovDI_with_multiplers}, we transform it into a DCS via Stewart's reformulation \cite{Stewart1990b}. In this reformulation, it is assumed that the sets $R_i$ are represented via the \textit{discriminant functions} $g_i(\cdot)$: 
\begin{align}\label{eq:stewart_sets}
	R_i  = \{ x \in \R^{n_x} \mid  g_i(x) < \min_{j\in \I\!,\, j \neq i } g_j(x)\}.	
\end{align}
Given the more intuitive representation via the sign matrix $S$ in Eq. \eqref{eq:standard_sets}, it can be shown that the function $g: \R^{n_x}  \to \R^{\Nsys}$ whose components are $g_i(x)$ can be found as \cite{Nurkanovic2022}:
\begin{align}\label{eq:indicator_func_formula}
		g(x) = - S c(x).
\end{align}
With this representation, the convex multipliers in the r.h.s. of \eqref{eq:FilippovDI_with_multiplers} can be found as  a solution of a suitable Linear Program (LP) \cite{Stewart1990b}, and \eqref{eq:FilippovDI_with_multiplers} is equivalent to
\begin{subequations}
\begin{align}
\dot{x}= F(x,u) \theta(x),\\
\theta(x)  \in\arg\min_{\tilde{\theta} \in \R^{\Nsys}} \quad & g(x)^\top \, \tilde{\theta} 
	\quad	\textrm{s.t.} \quad  e^\top \tilde{\theta} = 1
	,\ \tilde{\theta}\geq0 \label{eq:stewart_lp}.
\end{align}
\end{subequations}
We use a C-function $\Phi(\cdot,\cdot)$ for the complementarity conditions and write the KKT conditions of the LP \eqref{eq:stewart_lp} as a nonsmooth equation
\begin{align}\label{eq:lp_kkt}	
	G_{\lp}(x,\theta,\lambda,\mu) & \coloneqq \begin{bmatrix}
		g(x) - \lambda - \mu e\\
		1- e^\top \theta\\
		\Phi( \theta,\lambda) 
	\end{bmatrix} = 0,
\end{align}
where $\lambda \in \R^{n_f}_{\geq0}$ and $\mu \in \R$ are the Lagrange multipliers associated with the constraints of the LP  \eqref{eq:stewart_lp}. Note that $\mu = \min_{j\in\I} g_j(x)$.
Finally, the Filippov system is equivalent to the following DCS, which can be interpreted as a nonsmooth differential algebraic equation:
	\begin{align}\label{eq:dcs_1}
		\dot{x} & = F(x,u)\theta,\quad
    			0 = G_{\lp}(x,\theta,\lambda,\mu).
	\end{align}
A fundamental property of the multipliers $\lambda(\cdot)$ and $\mu(\cdot)$ is their continuity in time \cite[Lemma 5]{Nurkanovic2022}, whereas $\theta(\cdot)$ is in general a discontinuous function in time. 
 
\section{The Standard and FESD Discretizations for a Single Control Interval}\label{sec:algoritmic_ing}
\color{black}
This section describes the discretization of a single control interval in \NOSNOC\ via standard RK methods and FESD. We start with a standard RK method for the DCS \eqref{eq:dcs_1}. We subsequently introduce step-by-step the additional constraints which lead to FESD.
\subsection{Standard Runge-Kutta Discretization}
We consider a single control interval $[0,T]$ with a given constant control input $q$ and a given initial value $x_0 = s_0$. We divide the control interval into $\NFE$ finite elements (i.e., integration intervals) $[t_{n},t_{n+1}]$ via the grid points $0= t_0 < t_1 < \ldots <t_{\NFE} = T$. On each of these intervals, we apply an $\Nstg$-stage RK scheme, which is defined by its Butcher tableau entries $a_{i,j},\; b_i,\; c_i,\; i,j\in \{1,\ldots,\Nstg\}$ \cite{Hairer1991}. We denote the step-size as $h_{n} = t_{n+1} - t_{n},\; n = 0, \ldots,\NFE-1$. The approximation of the state at the grid points $t_n$ is denoted by $x_n \approx x(t_n)$. The time derivative of the state at the stage points $t_n + c_i h_n,\; i = 1,\ldots, \Nstg$, for a single finite element are collected in the vector $V_n \coloneqq (v_{n,1}, \ldots, v_{n,\Nstg}) \in \R^{\Nstg \cdot n_x}$. The stage values for the  algebraic variable  $\theta(\cdot)$ are collected in $\Theta_n \coloneqq (\theta_{n,1}, \ldots, \theta_{n,\Nstg} )\in \R^{\Nstg \cdot \Nsys}$. The vectors 
$\Lambda_n \in \R^{\Nstg \cdot \Nsys}$ and $M_n \in \R^{\Nstg}$ are defined accordingly.
Let $x_n^{\mathrm{next}}$ denote the value at the next time step $t_{n+1}$, which is obtained after a single RK step.

Now we can write the RK equations for the DCS \eqref{eq:dcs_1} in a compact \textit{differential} form. We summarize all RK equations of a finite element in $G_{\irk}(x_n^{\mathrm{next}},Z_n,h_n,q)= 0$, where $Z_n =(x_n,\Theta_n,\Lambda_n,M_n,V_n)$ collects all internal variables, and define 
\begin{align*}
\begin{split}
&G_{\irk}(x_n^{\mathrm{next}},Z_n,h_n,q)\! \coloneqq\!\!\! \\
&\begin{bmatrix}
\! v_{n,1}\! -\!  F(x_n +h_n \sum_{j=1}^{\Nstg} a_{1,j} v_{n,j},q)\theta_{n,1}\\
\vdots\\
v_{n,\Nstg} \! -\! F(x_n +h_n \sum_{j=1}^{\Nstg} a_{\Nstg,j} v_{n,j},q)\theta_{n,\Nstg}\\
G_{\lp}(x_n + h_n\sum_{j=1}^{\Nstg} a_{1,j} v_{n,j},\theta_{n,1},\lambda_{n,1},\mu_{n,1})\\
\vdots\\
G_{\lp}(x_n + h_n\sum_{j=1}^{\Nstg} a_{\Nstg,j} v_{n,j},\theta_{n,\Nstg},\lambda_{n,\Nstg},\mu_{n,\Nstg})\\
x_n^{\mathrm{next}} - x_n - h_n \sum_{i=1}^{\Nstg} b_i v_{n,i}
	\end{bmatrix}.
\end{split}
\end{align*}
To summarize all conditions for a single control interval in a compact way, we introduce some new notation. The variables for all finite elements of a single control interval are collected in the following vectors $\mathbf{x}= (x_0,x_0^\mathrm{next},\ldots,x_{\NFE}) \in \R^{(2\NFE+1) n_x}$, $\mathbf{V} = (V_0,\ldots,V_{\NFE-1}) \in \R^{\NFE \Nstg n_x}$ and $\mathbf{h}\coloneqq (h_0,\ldots,h_{\NFE-1})\in \R^{\NFE}$. The vectors
$\mathbf{\Theta}\in \R^{\NFE\Nstg\Nsys}$, $\mathbf{\Lambda}\in \R^{\NFE\Nstg\Nsys}$ and $\mathbf{M}\in \R^{\NFE\Nstg}$ are defined analogously. The vector $\mathbf{Z} = (\textbf{x},\mathbf{\Theta},\mathbf{\Lambda},\mathbf{M},\mathbf{V})$ collects all \textit{internal} variables.

Finally, we can summarize all computations over a single control interval and interpret it as a discrete-time nonsmooth system:
\vspace{-0.4cm}
	\begin{align}\label{eq:dcs_irk}
		{s}_1 \! = \!F_{\mathrm{std}}(\mathbf{Z}),\,
		0\! =\! 	G_{\mathrm{std}}(\mathbf{Z},\mathbf{h},s_0,q)
	\end{align}
with $F_{\mathrm{std}}(\mathbf{Z})  = x_{\NFE}$ and 
\begin{align*}
G_{\mathrm{std}}(\mathbf{Z},\mathbf{h},s_0,q)
		\coloneqq	&
		\begin{bmatrix}
			x_0  - s_0\\
			G_{\irk}(x_0^{\mathrm{next}},Z_0,h_0,q)\\
			x_1- x_0^{\mathrm{next}}\\
			\vdots\\
			G_{\irk}(x_{\NFE-1}^{\mathrm{next}},Z_{\NFE-1},h_{\NFE-1},q)\\
			x_{\NFE}- x_{\NFE-1}^{\mathrm{next}}
		\end{bmatrix}.
\end{align*}
Note that we keep a dependency on ${h}_n$ in \eqref{eq:dcs_irk}, but ${h}_n$ is implicitly given by the chosen discretization grid. This also means that for a standard RK scheme for DCS, higher order accuracy can be achieved only if the grid points $t_{n}$ coincide with all switching points, which is in practice impossible to achieve.
\subsection{Cross-Complementarity}
In FESD, the step-sizes $h_n$ are left as degrees of freedom such that the grid points $t_n$ can coincide with the switching times. Consequently, the switches should not happen on the stages inside a finite element. To exploit the additional degrees of freedom and to achieve these two effects we introduce additional conditions to the RK equations \eqref{eq:dcs_irk} called \textit{cross complementaries}. A key assumption, of course, is that there are more grid points in the interior of the grid than switching points.

For ease of exposition, we focus on the case where the right-boundary point of a finite element is also an RK-stage point, i.e., $c_{\Nstg} = 1$ and $t_{n+1} = t_n + c_{\Nstg}h_n$. Extensions can be found in \cite{Nurkanovic2022}. To achieve implicit and exact switch detection at the boundaries of $[t_n,t_{n+1}]$ and to avoid switching inside an element we exploit the fact that $\lambda(\cdot)$ and $\mu(\cdot)$ are continuous functions. We need their values at $t_n$ and $t_{n+1}$ which are denoted by $\lambda_{n,0},\; \mu_{n,0}$ and $\lambda_{n,\Nstg},\; \mu_{n,\Nstg}$, respectively. Due to continuity, we impose that $\lambda_{n,\Nstg}=\lambda_{n+1,0}$ and $\mu_{n,\Nstg}=\mu_{n+1,0}$ and use only the right boundary points of the finite elements ($\lambda_{n,\Nstg}$ and $\mu_{n,\Nstg}$) in the sequel.
\color{black}
To achieve the effects described above, we introduce the cross complementarity conditions which read  as \cite{Nurkanovic2022}: 
\begin{align}\label{eq:cross_comp}
	&0\! = \! G_{\mathrm{cross}}(\mathbf{\Theta},\mathbf{\Lambda}) \! \! \coloneqq \! \!
	\begin{bmatrix}
		\sum_{i = 1}^{\Nstg}\! \sum_{j=1,j\neq i}^{\Nstg}\theta_{1,i}^\top \lambda_{1,j}\\
		\vdots\\
		\sum_{i = 1}^{\Nstg}\! \sum_{\substack{j=0,\\j\neq i}}^{\Nstg}\theta_{\NFE\!-\!1,i}^\top \lambda_{\NFE\!-\!1,j}
	\end{bmatrix}.
\end{align}
This additional constraint ensures two very important properties: 
(i) we have the same active-set in \eqref{eq:dcs_irk} in $\Phi(\theta_{n,m},\lambda_{n,m})$ for all $m$ and changes can happen only for different $n$, i.e., at grid points $t_{n}$, 
(ii) whenever the active-sets for two neighboring finite elements differ in the $i$-th and $j$-th components of $\Phi(\theta_{n,m},\lambda_{n,m})$, then these two components of $\lambda_{n,\Nstg}$ must be zero \cite{Nurkanovic2022}.
This will implicitly result in the constraint $0= g_i(x_{n+1})-g_j(x_{n+1})$ (which comes from \eqref{eq:lp_kkt} and the fact that $\mu_{n,\Nstg} = \min_j g_j(x_{n+1})$). This defines the boundary between regions and $R_i$ and $R_j$, cf. \eqref{eq:stewart_sets}. Thus, it implicitly forces $h_n$ to adapt for exact switch detection.
\subsection{Step-Equilibration}
If no switches occur then also no active-set changes happen, hence the constraints \eqref{eq:cross_comp} are trivially satisfied. Consequently, the step-size $h_n$ can vary in a possibly undesired way and the optimizer can play with the discretization accuracy. To remove the spurious degrees of freedom we introduce an indicator function $\eta(\cdot)$ evaluated at the \textit{inner} grid points $t_n,\; n = 1,\dots,\NFE-1$ and its value at $t_n$ is denoted by $\eta_n$. It has the following property: if a switch happens at $t_n$ its value is zero, otherwise it is strictly positive. We omit the details on how a function $\eta(\cdot)$ is derived and refer to \cite{Nurkanovic2022}.
The discrete-time function $\eta(\cdot)$ depends on the values of $\Theta_{n}$ and $\Lambda_{n}$ of neighboring finite elements and we define
\begin{align*}
	\eta_n (\mathbf{\Theta},\mathbf{\Lambda}) \coloneqq \eta(\Theta_{n-1},\Lambda_{n-1},\Theta_n,\Lambda_n).
\end{align*}
Thus, the constraint $0=G_{\mathrm{eq}}(\mathbf{h},\mathbf{\Theta},\mathbf{\Lambda})$ removes the possible spurious degrees of freedom in $h_n$, where:
\begin{align}\label{eq:step_eq}
&G_{\mathrm{eq}}(\mathbf{h},\mathbf{\Theta},\mathbf{\Lambda}) 
\!\!	\coloneqq\!\!\!
	\begin{bmatrix}
		(h_{1}-h_{0})\eta_1(\mathbf{\Theta},\!\mathbf{\Lambda}) \\
		\vdots\\
		(h_{\NFE\!-\!1}-h_{\NFE\!-\!2})\eta_{\NFE\!-\!1}(\mathbf{\Theta},\!\mathbf{\Lambda}) \\
	\end{bmatrix}\!\!.
\end{align} 
We call the condition \eqref{eq:step_eq} \text{step-equilibration}. A consequence of \eqref{eq:step_eq} are locally equidistant state discretization grids between switching point, within a single control interval.
Since this constraint can be quite nonlinear, \NOSNOC\ offers several reformulations and heuristics that help numerical convergence.
\subsection{Finite Elements with Switch Detection}\label{sec:fesd}
We now use the ingredients explained above to state the FESD method. Similar to the standard RK scheme \eqref{eq:dcs_irk}, we summarize all computations over a single control interval and interpret it as a discrete-time nonsmooth system where internally exact switch detection is happening. The next step is computed by
\begin{align}\label{eq:fesd_compact}
	s_{1} \!=\! F_{\fesd}(\mathbf{Z}),\ \!	0 \!= \!G_{\fesd}(\mathbf{Z},\mathbf{h},s_0,q,T),
\end{align}
and $F_{\fesd}(\mathbf{Z})\!=x_{\NFE}$ renders the state transition map and the equation $0=G_{\fesd}(\mathbf{x},\mathbf{Z},q)$ collects all other internal computations including all RK steps within the regarded control interval:
\begin{align*}
	&G_{\fesd}(\mathbf{Z},\mathbf{h},s_0,q,T)\! \coloneqq\!\!
	\begin{bmatrix}
	G_{\mathrm{std}}(\mathbf{Z},\mathbf{h},s_0,q)\\
	G_{\mathrm{cross}}(\mathbf{\Theta},\mathbf{\Lambda})\\
	G_{\mathrm{eq}}(\mathbf{h},\mathbf{\Theta},\mathbf{\Lambda})\\
	\sum_{n=0}^{\NFE-1} h_n - T			
	\end{bmatrix}.
\end{align*}
 
The last condition ensures that the length of the considered time-interval is unaltered. In contrast to \eqref{eq:dcs_irk}, $h_n$ are now degrees of freedom, $s_0,q$ and $T$ are given parameters.
\color{black}
 The formulation \eqref{eq:fesd_compact} can be used as an integrator with exact switch detection for PSS \eqref{eq:ocp_pss}. This feature is implemented in \NOSNOC\ via the function \texttt{integrator\_fesd()}. It can automatically handle all kinds of switching cases such as: crossing a discontinuity, sliding mode, leaving a sliding mode or spontaneous switches \cite{Filippov2013}. 

\section{Discretizing and Solving a Nonsmooth Optimal Control Problem}\label{sec:ocp_discretization}
This section outlines how a nonsmooth OCP is discretized in \NOSNOC\ and how the resulting MPCC is solved.
\subsection{Multiple Shooting-Type Discretization with FESD}
One of the main goals of \NOSNOC\ is to numerically solve a discretized version of the  OCP \eqref{eq:ocp}. We consider $\Nctrl\geq 1$ control intervals of equal length, indexed by $k$, with piecewise constant controls collected in $\mathbf{q} = (q_0,\ldots,q_{\Nctrl-1})\in \R^{\Nctrl n_u}$. All internal variables are additionally equipped with an index $k$. On every control interval $k$, we apply an FESD discretization \eqref{eq:fesd_compact} with $N_{\mathrm{FE}}$ internal finite elements. The state values at the control interval boundaries are collected in  $\mathbf{s} = (s_0,\ldots,s_{\Nctrl})\in\R^{(\Nctrl+1)n_x}$. 
The vector ${\mathcal{Z}} = (\mathbf{{Z}}_0,\ldots,\mathbf{{Z}}_{\Nctrl-1})$ collects all internal variables and $\mathcal{H} = (\mathbf{h}_0,\ldots,\mathbf{h}_{\Nctrl-1})$ all step-sizes.
Finally the discretized OCP reads as:
\begin{subequations}\label{eq:ocp_discrete_time}
	\begin{align}
		\min_{\mathbf{s},\mathbf{q},\mathcal{Z},\mathcal{H}} \quad & \sum_{k=0}^{\Nctrl-1} \hat{f}_q(s_k,\mathbf{x}_k,q_k)+  \hat{f}_{\mathrm{T}}(s_{\Nctrl}) \\
		\textrm{s.t.} \quad  &s_{0} = \bar{x}_0,\\
		&{s}_{k+1}  = F_{\fesd}(\mathbf{x}_k),\;  k = 0,\ldots,\Nctrl\!-\!1,\\
		&0 = G_{\fesd}(\mathbf{x}_k,\mathbf{Z}_k,q_k),\; \! k = 0,\ldots,\Nctrl\!\!-\!\!1,\\
		&0 \geq G_{\mathrm{ineq}}(s_k,q_k),\; k = 0,\ldots,\Nctrl-1,\\
		&0 \geq G_{\mathrm{T}}(s_{\Nctrl}),
	\end{align}
\end{subequations}
where $\hat{f}_q:\R^{n_x}\times \R^{(\NFE+1)\Nstg n_x} \times \R^{n_u}\to \R$ and $\hat{f}_{\mathrm{T}}:\R^{n_x} \to \R$ are the discrteized stage and terminal costs, respectively.
\subsection{Reformulating and Solving MPCC}\label{sec:mpcc}
The discrete-time OCP \eqref{eq:ocp_discrete_time} is an MPCC. It can be written more compactly as
\begin{subequations} \label{eq:mpcc_generic}
	\begin{align}
		\min_{w} \quad & f(w)\\
		\textrm{s.t.} \quad 0 & \leq h(w),\\
		0&\leq w_1 \perp w_2\geq 0 \label{eq:mpcc_generic_cc},
	\end{align}
\end{subequations} 
where $w=(w_0,w_1,w_2) \in \mathbb{R}^{n_w}$ is a given decomposition of the problem variables. MPCC are difficult nonsmooth NLP which violate e.g., the MFCQ at all feasible points \cite{Anitescu2007}. Fortunately, they can often be solved efficiently via reformulations and homotopy approaches \cite{Anitescu2007,Scholtes2001}. We briefly discuss the different ways of solving MPCC that are implemented in \NOSNOC. They differ in how Eq. \eqref{eq:mpcc_generic_cc} is handled. In all cases, $w_1,w_2\geq0$ is kept unaltered and the bilinear constraint $w_1^\top w_2 =0$ is treated differently.

In a homotopy procedure, we solve a sequence of more regular, relaxed NLP related to \eqref{eq:mpcc_generic} and parameterized by a homotopy parameter $\sigma_i \in \R_{\geq0}$. Every new NLP is initialized with the solution of the previous one.
 In all approaches the homotopy parameter is updated via the rule: $\sigma_{i+1} = \kappa \sigma_i,\, \kappa \in (0,1),\, \sigma_0 >0$,  where $i$ is the index of the NLP in the homotopy.
In the limit as $\sigma_i \to 0$ (or often even for a finite $i$ and $\sigma_i$) the solution of the relaxed NLP matches a solution of \eqref{eq:mpcc_generic}. \NOSNOC\ supports the following approaches:
\paragraph*{Smoothing and Relaxation} In smoothing the bilinear term is replaced by the simpler constraint $w_1^\top w_2 = \sigma_i$ and in \textit{relaxation} by $w_1^\top w_2 \leq \sigma_i$. Under certain assumptions for $\sigma_i \rightarrow 0$ a solution of the initial MPCC \eqref{eq:mpcc_generic} is obtained \cite{Scholtes2001}.
\paragraph*{$\ell_1$-Penalty} In this approach, the bilinear constraint is discarded and the term $\frac{1}{\sigma_i} w_1^\top w_2$ is added to the objective, which is a penalized $\ell_1$ norm of the complementarity residual. When the penalty $\frac{1}{\sigma_i}$ exceeds a certain (often finite) threshold we have $w_1^\top w_2=0$ and the solution of such an NLP is a solution to \eqref{eq:mpcc_generic} \cite{Anitescu2007}.
\paragraph*{Elastic Mode}
In \textit{elastic mode} (sometimes called $\ell_\infty$-approach) \cite{Anitescu2007}, a bounded scalar slack variable $\gamma \in [0,\bar{\gamma}]$ is introduced. The relaxed bilinear constraint reads as $w_1^\top w_2 \leq \gamma$ and we add to the objective $\frac{1}{\sigma_i}\gamma$. Variants with $w_1^\top w_2 = \gamma$ and $-\gamma \leq w_1^\top w_2 \leq \gamma$ are supported as well. Once the penalty $\frac{1}{\sigma_i}$ exceeds a certain (often finite) threshold, we have $\gamma =0$ and we recover a solution of \eqref{eq:mpcc_generic} \cite{Anitescu2007}.

\section{NOSNOC Tutorials and a Benchmark}\label{sec:nosnoc}
In this section, we provide two short tutorials on the use of \NOSNOC. A numerical benchmark where we compare our software to conventional approaches is presented as well.
\subsection{Solving a Time-Optimal Control Problem}
\begin{figure}[t]
	\centering
	\vspace{0.05cm}
	\includegraphics[scale=0.65]{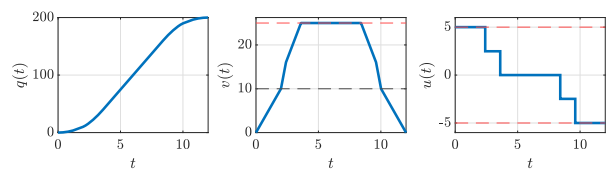}
	\caption{The position of the car $q(t)$ is shown in the left plot, the velocity $v(t)$ in the middle plot. Note the increase in  acceleration in the turbo mode for $v>\bar{v}$. The right plots shows the optimal control $u(t)$.}
	\label{fig:car}
\end{figure}

We regard a time-optimal control problem of a double-integrator car model with a normal and turbo mode. The state vector $x = (q,v)\in \R^2$ consists of the car's position $q$ and velocity $v$. 
The PSS reads as
\begin{align}\label{eq:pss_car}
	\dot{x} &= \begin{cases}
		(v,u),&\; \mathrm{if} \;v < \bar{v}\\
		(v,3u),&\; \mathrm{if} \;v >\bar{v}
	\end{cases}.
\end{align}
Following Section \ref{sec:prelimit}, we have $f_1(x,u) = (q,u)$ (nominal), $f_2(x,u) = (q,3u)$ (turbo). The two regions $R_1$ and $R_2$ described by $c(x) = v-\bar{v}$ and $S = \begin{bmatrix}
	-1 &  1
\end{bmatrix}^\top$.
The car should reach the state $x_{\mathrm{goal}} = (200,0)$ in optimal time $T$. Additionally, we have constraints on the velocity $|v| \leq v_{\mathrm{max}}$ and control $|u|\leq u_{\mathrm{max}}$. The parameters are $v_{\mathrm{max}} = 25$, $u_{\mathrm{max}} = 5$ and $\bar{v} = 10$. This OCP is formulated and solved with \NOSNOC\ using the code:
 \begin{lstlisting}
[settings] = default_settings_nosnoc();  
settings.time_optimal_problem = 1;
settings.n_s = 2; 
model.N_stg = 10; model.N_FE = 3; model.T = 1;    
q = MX.sym('q'); v = MX.sym('v'); 
model.x = [q;v]; model.x0 = [0;0];
model.lbx = [-inf;-25]; model.ubx = [inf;25];
u = SX.sym('u'); model.u = u;
model.lbu = -5; model.ubu = 5;
f_1 = [v;u]; f_2 = [v;3*u];
model.F = [f_1 f_2];  
model.c = v-10; model.S = [-1;1]; 
model.g_terminal = [q-200;v-0];
[results,stats,model,settings] = nosnoc_solver(model,settings);
\end{lstlisting}
The function \texttt{default\_settings\_nosnoc()} returns a MATLAB \texttt{struct} with default values for all possible tuning parameters. The needed time-transformations are automated by the flag \texttt{settings.time\_optimal\_problem = 1}. 
For the FESD-RK method we keep the default choice of a Radau II-A, hence we have with $\Nstg=2$ an accuracy order of 3~\cite{Hairer1991}. 
The MATLAB \texttt{struct} named \texttt{model} stores user input data, given in lines 4 to 13, which defines the OCP \eqref{eq:ocp}. 
\NOSNOC\ automates all definitions, reformulations and updates the \texttt{model} with all CasADi expressions for the DCS \eqref{eq:dcs_1}.  
Moreover, possible inconsistencies in the provided \texttt{settings} are refined.
Finally, in line 14  we solve the discretized OCP with a homotopy as described in Section \ref{sec:mpcc}. The solution trajectory is given in Fig. \ref{fig:car}. 
The user has access to all tuning parameters, intermediate results for all homotopy iterations and to all CasADi symbolic expressions and \texttt{Function} objects. This facilitates rapid prototyping and detailed analysis of solutions.
 
\begin{figure}[t]
	\centering
	\vspace{-0.2cm}
	\includegraphics[scale=0.51]{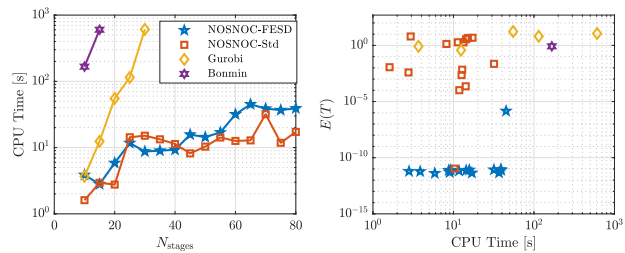}
	\vspace{-0.3cm}
	\caption{Comparison of \NOSNOC\ to mixed integer formulations. The left plot show CPU time as function of number of control intervals $\Nctrl$. The right plot show the solution accuracy as function of CPU time in a Pareto plot.}
	\label{fig:bencmark}
\end{figure}
\subsection{Numerical Benchmark}
We solve the OCP from the last section with four different approaches. We use \NOSNOC\ with the FESD discretization \eqref{eq:fesd_compact} and \NOSNOC\ with the standard discretization \eqref{eq:dcs_irk}. The latter approach is closely related to the smoothing approach in \cite{Stewart2010} and \cite{Nurkanovic2020}.  For the MPCC, we use in both cases the relaxation approach as it is usually the most robust one. 
Additionally, we make a \textit{big M} reformulation of the PSS \eqref{eq:pss_car} and solve a mixed integer nonlinear program (MINLP). Switches are allowed only at the control interval boundaries, hence we have two binary variables per control interval. We solve the MINLP with the dedicated solver Bonmin~\cite{Bonami2005}. Moreover, since the only non-linearity is in time $T$, we fix it and make a bisection-type search in $T$. For every fixed $T$ we solve a MILP with Gurobi. The MILP with the smallest $T$ that is still feasible delivers the optimal solution. We vary $\Nctrl$ from 10 to 80 with steps of 5. The computations are aborted if the time-limit of 10 minutes is exceeded.

The results of the benchmark are depicted in Fig. \ref{fig:bencmark}. \NOSNOC-FESD is slightly slower than \NOSNOC-Std, since it has $\Nctrl\NFE$ more variables, as $h_n$ are degrees of freedom. We compare also the solution quality by making a high-accuracy simulation $x_{\mathrm{sim}}(t)$ of \eqref{eq:pss_car} with the obtained optimal controls $\mathbf{q}$. We compare the terminal constraint satisfaction $E(T) = \| x_{\mathrm{sim}}(T) -x_{\mathrm{goal}}\|$. Due to the exact switch detection property, \NOSNOC-FESD has by far the most accurate solutions. The outlier where \NOSNOC-Std achieves high accuracy corresponds to a local minima without switches. We see that even a simple nonsmooth OCP is difficult to solve with conventional approaches, whereas  \NOSNOC-FESD provides faster and several orders of magnitude more accurate solutions. 
Further detailed comparisons of FESD to the standard approach can be found in \cite{Nurkanovic2022}. 
\subsection{An Example with State Jumps and Time-Freezing}\label{sec:time_freezing_example}
\begin{figure}[t]
	\centering
	\vspace{0.05cm}
	\includegraphics[scale=0.52]{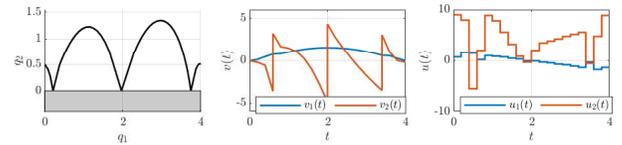}
	\vspace{-0.3cm}
	\caption{The illustration of the optimal solution $q(t)$ is shown in the left plot. The middle plot shows the optimal velocities $v(t)$ as a function of the \textit{physical time} $t$, where the state jumps are recovered. The right plot shows the optimal controls $u(t)$.}
	\label{fig:ocp_time_freezing_sol}
\end{figure}

In this subsection we illustrate how to use \NOSNOC\ with systems with state jumps. We consider a planar bouncing ball with elastic impacts.  The state vector is defined as $x = (q,v) \in \R^4$ with $q = (q_1,q_2)\in \R^2$ and $ v= (v_1,v_2)\in \R^2$ being the ball's position and velocity, respectively. The initial state is $x(0) = (0,0.5,0,0)$ and the ball is controlled with some force $u \in \R^2$. The ODE with state jumps reads as:
\begin{subequations}\label{eq:ode_state_jump}
\begin{align}
	\dot{q} &= v,\; \dot{v} = u - (0,g),\\
	v_2(t^+) &= - ev_2(t^-),\; \text{if}\; q_2(t) = 0\ \text{and}\ v_2(t)<0.
\end{align}
\end{subequations}
where $e \in \left(0,1\right]$ is the coefficient of restitution and determines the post impact velocity. The goal is to reach $q_{\mathrm{f}}= (4,0.5)$ with a minimal quadratic control effort modeled with the stage cost $f_q(x,u) =u^\top u$ and a minimal terminal velocity expressed via $f_{\mathrm{T}}(x) = 100v^\top v$, with $T=4$. The control force is bounded such that it is weaker than the gravitational force, i.e., $u^\top u \leq u_{\mathrm{max}}^2$. The chosen parameters are $e= 0.9$, $g =9.81$, $u_{\mathrm{max}} =9$. The following \NOSNOC\ code solves the described nonsmooth OCP with state jump:
 \begin{lstlisting}
[settings] = default_settings_fesd(); 
settings.time_freezing = 1; settings.n_s = 3; 
model.T = 4; model.N_stg = 20; model.N_FE = 3;
q = MX.sym('q',2); v = MX.sym('v',2);
model.x = [q;v]; model.x0 = [0;0.5;0;0];
u = MX.sym('u',2);  model.u = u;
model.c = q(2);  model.e = 0.9;
model.f = [v;u-[0;9.81]]; 
model.f_q = u'*u; model.f_q_T = 100*v'*v;
model.g_ineq = u'*u-9^2;
model.g_terminal = q-[4;0.5];
[results,stats,model,settings] = nosnoc_solver(model,settings);
\end{lstlisting}
The flag \texttt{settings.time\_freezing = 1} ensures that system with state jumps \eqref{eq:ode_state_jump} is transformed into a PSS of the form of \eqref{eq:ocp_pss} via the time-freezing reformulation \cite{Nurkanovic2021}. A solution trajectory is given in Figure \ref{fig:ocp_time_freezing_sol}, note the state jumps in $v_2(t)$ in the middle plot.

Many more settings can be changed by the user, for example, one can choose between different MPCC reformulations via \texttt{mpcc\_mode}, control the sparsity of the cross complementarities \texttt{cross\_complementarity\_mode} and so on. A few more examples and a detailed user manual are available \NOSNOC's repository \cite{Nurkanovic2022c}. 
\color{black}
\section{Conclusion and Outlook}\label{sec:conclusions}
In this letter we presented \NOSNOC, an open-source software package for nonsmooth numerical optimal control. With the help of the Finite Elements with Switch Detection (FESD) method and the time-freezing reformulation, it enables practical and high accuracy optimal control of several different classes of nonsmooth system in a unified way. The discretized OCP are solved with techniques solely from continuous optimization, without the need for any integer variables.  
All reformulations and details are hidden but accessible such that a convenient use for users with different knowledge levels of the field is ensured.
\color{black}

In future work, we aim to implement a python version of \NOSNOC. Moreover, further algorithmic developments in FESD, e.g., different reformulations of PSS into DCS, support for time-freezing for other classes of hybrid systems will be implemented.
\bibliographystyle{ieeetran}
\bibliography{bib/syscop1}
\end{document}